\documentclass[11pt]{amsart}

\usepackage{amsmath,amsthm,amsfonts,amssymb,latexsym,amscd,pictexwd,dcpic}

\usepackage{eucal}
\usepackage{mathrsfs}
\usepackage{amssymb}

\usepackage{setspace}

\usepackage[paper=letterpaper,left=1.5in,right=1in,bottom=1in,top=1in,headheight=12pt,headsep=24pt,includefoot=true,includehead=true]{geometry}

\title{Topological Hochschild Homology of $K/p$ as a $K_p^\wedge$ module}
\author{Samik Basu$\sp 1$}
\address{Department of mathematics, University of Copenhagen,Denmark}
\email{basu@math.ku.dk}

\numberwithin{equation}{section}

\theoremstyle{definition}
\newtheorem{thm}{Theorem}[section]

\newtheorem{prop}[thm]{Proposition}
\theoremstyle{definition}
\newtheorem{defn}[thm]{Definition}
\newtheorem{rem}[thm]{Remark}
\numberwithin{equation}{section}
\newtheorem{ex}[thm]{Example}

\newcommand{\Maps}{\mathit{Maps}}
\newcommand{\End}{\mathit{End}}
\newcommand{\Aut}{\mathit{Aut}}
\newcommand{\units}{\mathit{units}}

\newcommand{\mode}{\mathit{mod}}
\newcommand{\THH}{\mathit{THH}}
\newcommand{\unit}{\mathit{unit}}

\begin{document}

\begin{abstract}

Let $R$ be an $E_\infty$-ring spectrum. Given a map $\zeta$ from a space $X$ to $BGL_1R$, one can construct a Thom spectrum, $X^\zeta$, which generalises the classical notion of Thom spectrum for spherical fibrations in the case $R=S^0$, the sphere spectrum. If $X$ is a loop space ($\simeq \Omega Y$) and $\zeta$ is homotopy equivalent to $\Omega f$ for a map $f$ from $Y$ to $B^2GL_1R$, then the Thom spectrum has an $A_\infty$-ring structure. The Topological Hochschild Homology of these $A_\infty$-ring spectra is equivalent to the Thom spectrum of a map out of the free loop space of $Y$. 

This paper considers the case $X=S^1$, $R=K_p^\wedge$, the p-adic $K$-theory spectrum, and $\zeta = 1-p \in \pi_1BGL_1K_p^\wedge$. The associated Thom spectrum $(S^1)^\zeta$ is equivalent to the mod p $K$-theory spectrum $K/p$. The map $\zeta$ is homotopy equivalent to a loop map, so the Thom spectrum has an $A_\infty$-ring structure. I will compute $\pi_*THH^{K_p^\wedge}(K/p)$ using its description as a Thom spectrum.
\end{abstract}
\footnotetext[1]{supported in part by ERC-Adg TMSS-228082}
\maketitle
\tableofcontents

\parskip=2pt
\parindent=0.5cm

\markboth{}{}
\thispagestyle{myheadings}

\pagenumbering{arabic}

\section{Introduction}

 The goal of this paper is to use generalised Thom spectra to calculate the Topological Hochschild Homology of $K/p$ in the category of modules over $K_p^\wedge$. 

 Let $R$ be a ring spectrum and $GL_1R$ its space of units. It is the $H$-space of homotopy automorphisms of $R$ as an $R$-module. An $R$-twisting of a space $X$ is a continuous map $\zeta$ from $X$ to $BGL_1R$. Associated to $\zeta$, one can define the Thom spectrum of $\zeta$, $X^\zeta$ (see \cite{A-B-G-H-R}). This notion specialises for $R=S^0$ to the Thom spectrum of a spherical fibration. The homotopy groups of $X^\zeta$ is the group of twisted $R$ homology classes with respect to the twisting $\zeta$. 

 Suppose that $R$ is an $E_\infty$-ring spectrum. Then its space of units is an infinite loop space. Given a map $f:BG \rightarrow B^2GL_1R$, let $\zeta\simeq \Omega f:G\rightarrow BGL_1R$. Then the Thom spectrum $G^\zeta$ admits an $A_\infty$ $R$-algebra structure.

\subsection{$K/p$ as a module over $K_p^\wedge$} 

 Suppose that  $R=K_p^\wedge$, the spectrum of $p$-adic K-theory. Let $G$ be the group $S^1$. A twisting on $S^1$ is a map $\zeta: S^1 \rightarrow BGL_1K_p^\wedge$. This is classified by the group $\pi_1(BGL_1K_p^\wedge)\cong \pi_0(GL_1K_p^\wedge)\cong Z_p^\times$. If we choose $\zeta=1-p \in Z_p^\times$, then the Thom spectrum $(S^1)^\zeta \simeq K/p$, the $\mode$ $p$ $K$-theory spectrum. Moreover, the twisting $\zeta$ can be realised as a loop map, and so, for every way of writing $\zeta \simeq \Omega f$ we get an $A_\infty$-ring structure on $K/p$ as an $K_p^\wedge$-module. 

\subsection{Topological Hochschild Homology of Thom spectra}

Given a map $f$ from $X$ to $B^2GL_1R$, let $G \simeq \Omega X$ and $\zeta \simeq \Omega f : G \simeq \Omega X \rightarrow BGL_1R$. In this case, the Thom spectrum $G^\zeta$ has an $A_\infty$-ring structure. We write $\eta^*f$ for the composite 

\begindc{0}[5]
\obj(-15,15){}
\obj(0,-5){}
\obj(0,10){$LX$}
\obj(10,10){$LB^2GL_1R$}
\obj(30,10){$B^2GL_1R\times BGL_1R$}
\obj(30,0){$BGL_1R\times BGL_1R$}
\obj(50,0){$BGL_1R$}
\mor{$LX$}{$LB^2GL_1R$}{}
\mor{$LB^2GL_1R$}{$B^2GL_1R\times BGL_1R$}{$\cong$}{}
\mor{$B^2GL_1R\times BGL_1R$}{$BGL_1R\times BGL_1R$}{$\eta \times id$}
\mor{$BGL_1R\times BGL_1R$}(47,0){}
\enddc

where $\eta: \Sigma R \rightarrow R$ is induced from $S^1\stackrel{\eta}{\rightarrow}S^0$ via $S^1 \wedge R \rightarrow S^0 \wedge R \simeq R$. In the above situation, $\THH^R(G^\zeta) \simeq LX^{\eta^* f}$. The case $R=S^0$ was proved in \cite{B-C-S}. The same argument applies for general $R$.

Using this identification of $\THH$ as a Thom spectrum, we compute the Topological Hochschild Homology of $K/p$. For odd primes p, 
$$ \pi_*(\THH^{K_p^\wedge}(K/p))= \left\{
\begin{array}{rl}
(Z/(p^\infty))^i & \text{if } * =2k\\
0 & \text{if } * =2k+1 
\end{array} \right.$$ 
where $i$ is an integer between $1$ and $p-1$ depending on the choice of $f$ with $\zeta \simeq \Omega f$. 

Similar results were obtained before by Angeltveit in \cite{A1}. He used the B{\"o}kstedt spectral sequence (see \cite{E-K-M-M}, chapter IX).

We can also form $\mode$ $p$ $K$-theory as a Thom spectrum by starting with $X=S^3$, $R=K_p^\wedge$ and $\zeta= p \in \pi_3(BGL_1K_p^\wedge)= \pi_2(GL_1K_p^\wedge)= Z_p$. Again, this $\zeta$ can be realised as a loop map and we can compute $\THH$ of these $A_\infty$-ring structures in an analogous way. This gives the same results.

\section{The Thom spectrum}

The notion of a generalised Thom spectrum used here is discussed in detail in \cite{A-B-G-H-R}. The construction resembles a twisted version of the group ring. Given an extension of a group $G$ by the units in a field $k$,
$$(\tau): 1\rightarrow k^* \rightarrow E \rightarrow G \rightarrow 1$$
the algebra  $k^\tau[G]= \mathbb{Z}[E]\otimes_{\mathbb{Z}[k^*]} k$ is a twisted group ring. If the extension $\tau$ is trivial, one gets the group ring $k[G]$. Imitating this definition of a twisted group ring for spectra leads to the construction of the Thom spectrum. One replaces the field $k$ by an $E_\infty$-ring spectrum $R$, and the units $k^*$ by the space of units $GL_1R$ acting on $R$.   

\subsection{The space of units and the Thom spectrum}

The space of units of a ring spectrum is a generalisation of the group of units of a commutative ring, the set of invertible elements under multiplication. It is defined to be the components of $\Omega^\infty R$ that lie over the units in $\pi_0 (R)$. Following \cite{A-B-G-H-R}, we make the definition 

\defn
 Let $R$ be an $E_\infty$-ring spectrum. Its space of units $GL_1R$ is defined to be the homotopy pullback, 
 
\begindc{0}[5]
\obj(0,25){}
\obj(0,5){}
\obj(30,20){$GL_1R$} 
\obj(50,20){$\Omega^\infty(R)$}
\obj(30,10){$\pi_0(R)^\times$}
\obj(50,10){$\pi_0(R)$}
\mor{$GL_1R$}{$\Omega^\infty(R)$}{}
\mor{$GL_1R$}{$\pi_0(R)^\times$}{}
\mor{$\Omega^\infty(R)$}{$\pi_0(R)$}{}
\mor{$\pi_0(R)^\times$}{$\pi_0(R)$}{}
\enddc

It follows from the definition that the homotopy classes of maps from a space $X$ to $GL_1R$ are given by 
$$[X,GL_1R]= R^0(X)^\times $$ 
the units of the cohomology ring $R^0(X)= [X,\Omega^\infty R]$. 
 
From the pullback diagram one can read off the homotopy groups of $GL_1R$, 
$$ \pi_n(GL_1R) = \left\{
\begin{array}{rl}
\pi_n(R) & \text{if } n>0\\
\pi_0(R)^\times & \text{if } n=0
\end{array} \right.$$                

We note that $GL_1R$ is an $H$-space for any ring spectrum $R$. If $R$ is $E_\infty$, then $GL_1R$ is an infinite loop space: there is a connective spectrum $gl_1R$ with $0^{th}$-space is $GL_1R$ (Theorem 3.2 in \cite{A-B-G-H-R}). 

We can view $\Omega^\infty R$ as the space of endomorphisms $\End_R(R,R)$, in the topological  category of $R$-modules, and $GL_1R= \Aut_R(R,R) \subset \End_R(R,R)$ as the subset of weak equivalences. Therefore, the units $GL_1R$ is the space of homotopy automorphisms of $R$ in the category of $R$-modules. In this way, the infinite loop space $GL_1R$ acts on the spectrum $R$ by weak equivalences, and $R$ is a module over the $E_\infty$ ring spectrum $\Sigma^\infty GL_1R_+$.

\defn Given a map $\zeta : X \rightarrow BGL_1R$, let $P$ be the $GL_1R$ bundle classified by $\zeta$ described as the pullback,
 
 \begindc{0}[5]
\obj(0,25){}
\obj(-15,10){}
\obj(10,20){$P$} 
\obj(30,20){$EGL_1(R)$}
\obj(10,10){$X$}
\obj(30,10){$BGL_1(R)$}
\mor(10,20)(25,20){}
\mor(10,20)(10,10){}
\mor(30,20)(30,10){}
\mor(10,10)(25,10){$\zeta$}
\enddc
 and define the associated Thom spectrum to be  
 $$X^\zeta = \Sigma^\infty P_+\wedge{^L}_{\Sigma^\infty GL_1(R)_+}R$$

 In the above $\wedge^L$ denotes the derived smash product in the category of modules over the $E_\infty$-ring spectrum $\Sigma^\infty GL_1R_+$ as in \cite{E-K-M-M}. We note from section 7 of \cite{A-B-G-H-R}, that the Thom spectrum functor commutes with homotopy colimits, and from section 8.6 of \cite{A-B-G-H-R} that it generalises the classical Thom spectrum of a spherical fibration.  

The Thom spectrum of the map $* \rightarrow BGL_1R$ is weakly equivalent to $R$, since the universal bundle associated to the inclusion of a point in $BGL_1R$ is isomorphic to $GL_1R$ and $\Sigma^\infty GL_1R_+ \wedge^L_{\Sigma^\infty GL_1R_+} R \simeq R$. 

Similarly, the Thom spectrum of a map $X\rightarrow BGL_1R$ which is null homotopic is weakly equivalent to $R\wedge X_+$. Indeed, the universal bundle associated to the constant map is $X\times GL_1R$. Then the Thom spectrum is $\Sigma^\infty (X\times GL_1R)_+\wedge^L_{\Sigma^\infty GL_1R_+} R \simeq (\Sigma^\infty X_+ \wedge \Sigma^\infty GL_1R_+) \wedge^L_{\Sigma^\infty GL_1R_+} R \simeq R\wedge X_+$. 

Suppose that the space $X\simeq \Sigma Y$, the reduced suspension on $Y$. Then, a map $X\stackrel{\zeta}{\rightarrow} BGL_1R$ is described by a map $Y\stackrel{\widehat{\zeta}}{\rightarrow} GL_1R$, via $[X,BGL_1R]\cong [\Sigma Y,BGL_1R]\cong [Y,GL_1R]$. Such a $\widehat{\zeta}$ is a unit in $R^0(Y)$ which induces $u_\zeta : R\wedge Y_+ \rightarrow R$.

\prop Suppose that $\zeta$ is a map from $X\simeq\Sigma Y$ to $BGL_1R$. Then, the Thom spectrum $X^\zeta$ is equivalent to the homotopy colimit of $(R\leftarrow R\wedge Y_+ \rightarrow R)$ where one of the maps is the projection $p_Y$ and the other is $u_\zeta$.

\begin{proof}

The space $X$ is the homotopy colimit of $*\leftarrow Y\rightarrow *$, and this gives a homotopy pushout square of Thom spectra, 

\begindc{\commdiag}[5]
\obj(0,25){}
\obj(-15,5){}
 \obj(10,20){$Y^\zeta$}
\obj(30,20){$*^\zeta$}
\obj(10,10){$*^\zeta$}
\obj(30,10){$(\Sigma Y)^\zeta$}
\mor(10,20)(10,10){}
\mor(12,20)(29,20){}
\mor(12,10)(28,10){}
\mor(30,20)(30,10){}
\enddc
The Thom spectrum $*^\zeta$ is weakly equivalent to $R$ and $Y^\zeta \simeq R \wedge Y_+$, so the homotopy pushout can be written as

\begindc{\commdiag}[5]
\obj(0,25){}
\obj(-15,5){}
 \obj(10,20){$R \wedge Y_+$}
\obj(30,20){$R$}
\obj(10,10){$R$}
\obj(30,10){$(\Sigma Y)^\zeta$}
\mor(10,20)(10,10){}
\mor(12,20)(30,20){}
\mor(10,10)(28,10){}
\mor(30,20)(30,10){}
\enddc

From this, one obtains a  Mayer Vietoris sequence for calculating the homotopy groups 
$$\ldots \rightarrow \pi_*(R \wedge Y_+) \rightarrow \pi_*(R)\oplus \pi_* (R) \rightarrow \pi_*((\Sigma Y)^\zeta) \ldots$$
To compute the maps in this sequence, one must examine the $GL_1R$-bundle over $X\simeq \Sigma Y$. This restricts to trivial bundles over the two copies of the cone of $Y$ inside $X$ and on their intersection $Y$, the bundles are identified via the map $\widehat{\zeta}:Y \rightarrow GL_1R$. 

In the long exact sequence, there are two maps $R_*(Y_+) \rightarrow \pi_*(R)$. One of these maps is given by the map from $Y$ to a point($p_Y$) and the other is the map $u_\zeta$ defined in the preceeding paragraph.   
\end{proof}

\rem The proposition describes the homotopy groups of the Thom spectrum as twisted $R$-homology groups. An $R$-twisting on a space $X$ can be defined as a 1-cocycle in the sheaf (of groupoids) - $\lbrace \units\, \mathit{in}\, R^0(X) \rbrace$. The groupoid of units in $R^0$ is classified by the units $GL_1R$, and therefore, 1-cocycles on $X$ are equivalent to $[X,BGL_1R]$. Therefore, a twisting is given by a continuous map $\zeta$ from $X$ to $BGL_1R$.   

For $X=\bigcup U_i$ a 1-cocycle defines units over $U_i\cap U_j$ satisfying a cocycle condition on further intersections. A twisted $R$ homology class is an element in each $R_*(U_i)$, two of which are identified using the values of the 1-cocycle on the intersections. The abelian group of these classes is defined to be the twisted $R$-homology of $X$ with respect to the twisting $\zeta$. This is isomorphic to the homotopy groups of the Thom spectrum $X^\zeta$. The proposition above verifies this in the case $X=\Sigma Y$, where $X$ is the union of two contractible open sets.

\subsection{Computations of some Thom spectra}
 
\prop Suppose that $\zeta:S^1 \rightarrow BGL_1K_p^\wedge$ represents $1-p \in \pi_1(BGL_1(K_p^\wedge))= \pi_0(GL_1(K_p^\wedge))= Z_p^\times$. Then, $(S^1)^\zeta \simeq K/p$. 
 \begin{proof}
 By Proposition 2.3 with $Y=S^0$, the Thom spectrum is a homotopy pushout 

\begindc{\commdiag}[5]
\obj(0,25){}
\obj(-15,5){}
 \obj(10,20){$K_p^\wedge\vee K_p^\wedge$}
\obj(30,20){$K_p^\wedge$}
\obj(10,10){$K_p^\wedge$}
\obj(30,10){$(S^1)^\zeta$}
\mor(10,20)(10,10){}
\mor(13,20)(29,20){}
\mor(10,10)(29,10){}
\mor(30,20)(30,10){}
\enddc

Therefore, there is a cofibre sequence 
$$ K_p^\wedge \vee K_p^\wedge \rightarrow K_p^\wedge \vee K_p^\wedge \rightarrow (S^1)^\zeta$$
Proposition 2.3 also identifies the left map in the sequence in suitable coordinates, to be given by the matrix 
\[ \left( \begin{array}{ccc}
1 & 1 \\
1 & 1-p  \end{array} \right)\] 

Therefore, the cofibre sequence can be rewritten as  
$$K_p^\wedge \stackrel{p}{\rightarrow} K_p^\wedge \rightarrow (S^1)^\zeta$$
so that $(S^1)^\zeta \simeq {K_p^\wedge}/p \simeq K/p$.
\end{proof}
 
\rem Consider the map $\zeta:S^1 \rightarrow BGL_1((S^0)_p^\wedge)$ given by $(1-p)$ as in the previous proposition. Then,  $(S^1)^\zeta \simeq (S^0)_p^\wedge /p \simeq M_p$ is the $\mode$ $p$ Moore spectrum. In fact, for any $\zeta: S^1 \rightarrow BGL_1R$, $(S^1)^\zeta \simeq \mathit{cofibre}(1-\zeta : R \rightarrow R)$. This follows from the argument above.

\prop Let $\zeta: S^3 \rightarrow BGL_1K_p^\wedge$ represent the element $p$ of 
$$[S^3,BGL_1(K_p^\wedge)] = \pi_3(BGL_1(K_p^\wedge))= \pi_2(GL_1(K_p^\wedge))= \pi_2(K_p^\wedge)\cong Z_p$$
 Then $(S^3)^\zeta \simeq K/p$.

\begin{proof}

The space $S^3$ is homotopy equivalent to the suspension of $S^2$. Proposition 2.3 implies the homotopy pushout 

\begindc{\commdiag}[5]
\obj(0,22){}
\obj(-15,8){}
 \obj(10,20){$ K_p^\wedge \wedge S^2_+$}
\obj(30,20){$K_p^\wedge$}
\obj(10,10){$K_p^\wedge$}
\obj(30,10){$(S^3)^\zeta$}
\mor(10,20)(10,10){}
\mor(13,20)(30,20){}
\mor(10,10)(29,10){}
\mor(30,20)(30,10){}
\enddc
and the associated Mayer Vietoris cofibre sequence 
$$ K_p^\wedge \wedge (S^2) \vee K_p^\wedge \rightarrow K_p^\wedge \vee K_p^\wedge \rightarrow (S^3)^\zeta.$$
In suitable coordinates, the map in the Mayer Vietoris sequence is given by the matrix 
\[ \left( \begin{array}{ccc}
1 & 0 \\
1 & p  \end{array} \right)\] 
and the sequence can be rewritten as
$$\Sigma^2 K_p^\wedge \stackrel{p}{\rightarrow} K_p^\wedge \rightarrow (S^3)^\zeta $$
By Bott periodicity $\Sigma^2 K_p^\wedge \simeq K_p^\wedge$ so that $(S^3)^\zeta \simeq K_p^\wedge / p$, as claimed. 
\end{proof}

\subsection{Ring Structures}

Suppose $R$ is an $E_\infty$-ring spectrum so that $GL_1R$ is an infinite loop space. Given $f:X\rightarrow B^2GL_1R$, and $\zeta:G\simeq \Omega X\stackrel{\Omega f}{\rightarrow} BGL_1R$, the Thom spectrum $G^\zeta$ has an $A_\infty$-ring structure. This follows from \cite{B1} and \cite{A-B-G-H-R}. This raises the question when a map 
$$\zeta:G \rightarrow BGL_1R$$
from a monoid $G$ is homotopy equivalent to a loop map, i.e. $\zeta \simeq \Omega f$ for 
$$f:BG \rightarrow B^2GL_1R.$$
We have the standard maps 
$$\Sigma G \stackrel{\sigma}{\rightarrow} BG,\, \Sigma GL_1R \stackrel{\sigma}{\rightarrow} BGL_1R$$
so the question is if 
$$\sigma \circ \Sigma \zeta : \Sigma G \rightarrow B^2GL_1R$$
extends over $BG$, 

\begindc{\commdiag}[5]
\obj(0,25){}
\obj(-15,5){}
\obj(10,20){$\Sigma G$}
\obj(30,20){$\Sigma BGL_1(R)$}
\obj(10,10){$BG$}
\obj(30,10){$B^2GL_1(R)$}
\mor{$\Sigma G$}{$\Sigma BGL_1(R)$}{$\Sigma \zeta$}
\mor{$\Sigma BGL_1(R)$}{$B^2GL_1(R)$}{$\sigma$}
\mor{$BG$}{$B^2GL_1(R)$}{$f$}
\mor{$\Sigma G$}{$BG$}{$\sigma$}
\enddc
 
\prop Let $G=S^1$, $R=K_p^\wedge$ and $\zeta = 1-p$ as in Proposition 2.5, then $(S^1)^\zeta \simeq K/p$ has an $A_\infty$-ring structure. 

\begin{proof}
The classifying space of $S^1$ is $CP^\infty$ so, in this case the diagram above is

\begindc{\commdiag}[5]
\obj(-23,13){}
\obj(0,-1){}
\obj(0,10){$S^2$}
\obj(20,10){$\Sigma BGL_1(K_p^\wedge)$}
\obj(0,0){$CP^\infty$}
\obj(20,0){$B^2GL_1(K_p^\wedge)$}
\mor(0,10)(15,10){$\Sigma (1-p)$}
\mor(20,10)(20,0){$\sigma$}
\mor(1,0)(15,0){$f$}
\mor(0,10)(0,0){$\sigma$}
\enddc 
 The space $CP^\infty$ has a CW structure made of even dimensional cells so that all the cells are attached along odd dimensional spheres. The spectrum $K_p^\wedge$ has non trivial homotopy groups only in even dimensions and hence, so does $B^2GL_1K_p^\wedge$. Thus, all the obstructions to extending the map $\Sigma 1-p$ must vanish, which implies that there is an $A_\infty$-ring structure on the Thom spectrum $K/p$. 
\end{proof}

\prop Suppose that $G=S^3$, $R=K_p^\wedge$, and $\zeta =p$ as in Proposition 2.7, then the Thom spectrum has an $A_\infty$-ring structure. 
\begin{proof}
 The classifying space of $S^3$ is the infinite quarternionic projective space $HP^\infty$, and $\Sigma S^3 = S^4 \rightarrow BS^3= HP^\infty$ is obtained by attaching even cells along maps of odd dimensional spheres. Therefore the extension problem can always be solved.
\end{proof}

\section{Topological Hochschild Homology of Thom spectra}
 
In the last section, we observed that the Thom spectrum of a loop map carries an induced $A_\infty$ structure. In this setting, there is a convenient description of the Topological Hochschild Homology as a Thom spectrum in the ideas of \cite{B-C-S} and \cite{Sch}. In the following $G$ will be a group, $X$ a space, and $G$ homotopy equivalent to $\Omega X$ as $A_\infty$-spaces. $R$ will be an $E_\infty$ ring spectrum.
 
 The Thom spectrum of a map $G\rightarrow BGL_1R$ is a twisted $R$-module generated by $G$. If this is a loop map, the construction is that of a twisted group ring. Recall that the Hochschild Homology of group rings over a field is given by 
$$HH_*(k[G])\cong k \otimes H_*(G,G)$$ 
where $G$ acts on itself by conjugation. This is the homology of the Borel construction $G_{hG}\simeq EG\times_G G \simeq LBG$, the free loop space of $BG$, and so, $HH_*(k[G])\cong k\otimes H_*(LBG)$. The analogous statement for Topological Hochschild Homology is the classical result of B\"{o}kstedt and Waldhausen,
$$\THH(\Sigma^\infty \Omega X_+) \simeq \Sigma^\infty LX_+.$$ 
 
 In the category of $R$-modules, the theorem is $\THH^R( R\wedge \Omega X_+) \simeq R\wedge LX _+$, computing the Topological Hochschild Homology of the Thom spectrum of the constant map. More generally, let $f:X\rightarrow BGL_1R$ and $\zeta \simeq \Omega f :G \rightarrow BGL_1R$, the Thom spectrum has an $A_\infty$-ring structure, and the Topological Hochschild Homology is the Thom spectrum of a map from $LX$ to  $BGL_1R$. 

 In the second part of the section, we apply the theorem for $R=K_p^\wedge$ and $G=S^1$, in the computation of the previous section. This implies that the Thom spectrum is homotopy equivalent to the cofibre of a certain map $K_p^\wedge \wedge CP^\infty_+ \rightarrow K_p^\wedge \wedge CP^\infty_+$. 
 
\subsection{ Identifying Topological Hochschild Homology as a Thom spectrum}
 
 Recall that the free loop space $LY$ fits into a fibration 
 $$\Omega Y \rightarrow LY \rightarrow Y$$ 
  If $Y$ is an $H$-space, then the fibration splits as $LY\simeq Y\times \Omega Y$. This is an equivalence of $H$-spaces if $Y$ is homotopy commutative. 
 
 Let $f$ be a map from $X$ to $B^2GL_1R$ and $\eta: B^2GL_1R \rightarrow \Omega B^2GL_1R$ be induced from the Hopf map by
\begin{align}
B^2GL_1R \simeq \Maps(S^2,B^4GL_1R) & \stackrel{\eta^*}{\rightarrow} \Maps(S^3, B^4GL_1R) \nonumber \\
                              &\simeq \Maps(S^1,\Omega^2B^4GL_1R) \nonumber \\
                              &\simeq \Maps(S^1,B^2GL_1R) \nonumber \\
                              &\simeq \Omega B^2GL_1R.  \nonumber
\end{align}

Let $L^\eta f$ be the map from $LX$ to $BGL_1R$ defined by the diagram 

 \begindc{0}[5]
 \obj(0,35){}
 \obj(-15,5){}
 \obj(0,30){$LX$}
\obj(15,30){$LB^2GL_1(R)$}
\obj(40,20){$\Omega B^2GL_1(R)$}
\obj(40,30){$B^2GL_1(R)\times \Omega B^2GL_1(R)$}
\obj(40,10){$BGL_1(R)$}
\mor{$LX$}{$LB^2GL_1(R)$}{$Lf$}
\mor{$LB^2GL_1(R)$}{$B^2GL_1(R)\times \Omega B^2GL_1(R)$}{$\simeq $}
\mor{$B^2GL_1(R)\times \Omega B^2GL_1(R)$}{$\Omega B^2GL_1(R)$}{$\eta \times id$}
\mor{$\Omega B^2GL_1(R)$}{$BGL_1(R)$}{$\simeq $}
\mor{$LX$}{$BGL_1(R)$}{$L^\eta f$}
\enddc  
The map $\eta \times id : B^2GL_1R \times \Omega B^2GL_1R \rightarrow \Omega B^2GL_1R $ is the product of the maps $\eta$ and $\mathit{id}$ using the $H$-space structure of $\Omega B^2GL_1R$. Without proof, we state:

\thm There is a homotopy equivalence
$$\THH^R(G^\zeta)\simeq (LX)^{L^\eta f} $$

This was proved in the case of the sphere spectrum in \cite{B-C-S}, \cite{Sch}. A similar argument applies for any $E_\infty$-ring spectrum $R$. This will be accomplished in a future publication.

\subsection{The example of $G=S^1$ and $R=K_p^\wedge$} 
By Proposition 2.8, we have the commutative diagram,

\begindc{\commdiag}[5]
\obj(-15,15){}
\obj(0,-5){}
\obj(0,10){$S^2$}
\obj(20,10){$\Sigma BGL_1(K_p^\wedge)$}
\obj(0,0){$CP^\infty$}
\obj(20,0){$B^2GL_1(K_p^\wedge)$}
\mor(0,10)(15,10){$\Sigma 1-p$}
\mor(20,10)(20,0){$\sigma$}
\mor(1,0)(15,0){$f$}
\mor(0,10)(0,0){$\sigma$}
\enddc 
and write $\THH^{K_p^\wedge}(K/p,f)$ for the Topological Hochschild Homology corresponding to this $A_\infty$-ring structure.
\prop $$\THH^{K_p^\wedge}(K/p,f)\simeq (LCP^\infty)^{\widehat{f}}$$
where $\widehat{f}$ is the composite, 
$$LCP^\infty \stackrel{Lf}{\rightarrow} LB^2GL_1K_p^\wedge \simeq B^2GL_1K_p^\wedge \times BGL_1K_p^\wedge \stackrel{p_2}{\rightarrow} BGL_1K_p^\wedge$$

\begin{proof}
By Theorem 3.1, $\THH^{K_p^\wedge}(K/p,f)\simeq (LCP^\infty)^{L^\eta f}$. Since $\pi_1(K_p^\wedge)=0$, $\eta =0$ in this case. Hence, the proposition.
\end{proof} 

The focus of the rest of the paper will be the calculation of $\pi_*((LCP^\infty)^{\widehat{f}})\cong \THH^{K_p^\wedge}(K/p,f)$. First of all we note that:

\prop There is a long exact sequence 
$${K_p^\wedge}_* CP^\infty \rightarrow {K_p^\wedge}_* CP^\infty \rightarrow \pi_*\THH^{K_p^\wedge}(K/p,f)\rightarrow {K_p^\wedge}_{*-1} CP^\infty \ldots $$

\begin{proof}
Note that $CP^\infty$ is an infinite loop space, and hence homotopy commutative, which implies that $LCP^\infty \simeq \Omega CP^\infty \times CP^\infty\simeq S^1\times CP^\infty$. The space $S^1$ is a union of two contractible open sets whose intersection is $S^0$, so, there is a homotopy pushout, 

\begindc{\commdiag}[5]
\obj(0,16){}
\obj(-25,0){}
\obj(0,13){$C P^\infty \sqcup C P^\infty$}
\obj(20,13){$C P^\infty$}
\obj(0,3){$C P^\infty$}
\obj(20,3){$LC P^\infty$}
\obj(40,8){$(*)$}
\mor(5,13)(19,13){}
\mor(0,13)(0,3){}
\mor(3,3)(18,3){}
\mor(20,13)(20,3){}
\enddc
and hence, a homotopy pushout square of Thom spectra
 
\begindc{\commdiag}[5]
\obj(0,15){}
\obj(-25,-5){}
\obj(0,10){$(CP^\infty \sqcup CP^\infty)^{\widehat{f}}$}
\obj(20,10){$(CP^\infty)^{\widehat{f}}$}
\obj(0,0){$(CP^\infty)^{\widehat{f}}$}
\obj(20,0){$(LCP^\infty)^{\widehat{f}}$}
\mor(6,10)(17,10){}
\mor(0,10)(0,0){}
\mor(3,0)(16,0){}
\mor(20,10)(20,0){}
\enddc

The two maps $CP^\infty \rightarrow LCP^\infty$ in (*) are the inclusion of constant loops, so, the two compositions $CP^\infty \rightarrow LCP^\infty \rightarrow LB^2GL_1K_p^\wedge \rightarrow BGL_1K_p^\wedge$ are nullhomotopic and the Thom spectra are $\simeq K_p^\wedge \wedge CP^\infty_+$. The map from $CP^\infty \sqcup CP^\infty$ to $BGL_1K_p^\wedge$ factors through $CP^\infty \rightarrow BGL_1K_p^\wedge$ so, the Thom spectrum $(CP^\infty \sqcup CP^\infty)^{f} \simeq K_p^\wedge \wedge CP^\infty_+ \vee K_p^\wedge \wedge CP^\infty_+$. Therefore, the pushout can be written as: 

\begindc{\commdiag}[5]
\obj(0,15){}
\obj(-25,-5){}
\obj(0,10){$K_p^\wedge \wedge CP^\infty_+ \vee K_p^\wedge \wedge CP^\infty_+$}
\obj(25,10){$K_p^\wedge \wedge CP^\infty_+$}
\obj(0,0){$K_p^\wedge \wedge CP^\infty_+$}
\obj(25,0){$(LCP^\infty)^{Lf}$}
\mor(10,10)(21,10){}
\mor(0,10)(0,0){}
\mor(3,0)(21,0){}
\mor(25,10)(25,0){}
\enddc

This gives a Mayer Vietoris sequence on homotopy groups,  
$$ \ldots \rightarrow{K_p^\wedge}_*(CP^\infty) \oplus {K_p^\wedge}_*(CP^\infty) \rightarrow {K_p^\wedge}_*(CP^\infty) \oplus {K_p^\wedge}_*(CP^\infty) \rightarrow \pi_*((LCP^\infty)^{Lf})\ldots$$
To simplify, one needs to understand the left hand map i.e., how $K_p^\wedge \wedge CP^\infty_+ \vee K_p^\infty \wedge CP^\infty_+$ maps to the two different copies of $K_p^\wedge \wedge CP^\infty_+$ in the pushout square. For that one needs to examine the structure of $P^{\widehat{f}}$, the $GL_1K_p^\wedge$-bundle over $S^1\times CP^\infty$ classified by $\widehat{f}$. 

 Following the pushout square (*), we see that $P^{\widehat{f}}$ is obtained by identifying two trivial bundles over $CP^\infty$ after restricting over $CP^\infty \sqcup CP^\infty$, via a map $u:CP^\infty \sqcup CP^\infty \rightarrow GL_1K_p^\wedge$. The adjoint of $u$ is the map $\widetilde{u}$ in the diagram, 

\begindc{\commdiag}[5]
\obj(0,15){}
\obj(-10,-3){}
\obj(0,10){$CP^\infty \sqcup CP^\infty$}
\obj(20,10){$CP^\infty \vee CP^\infty$}
\obj(40,10){$S^1 \times CP^\infty$}
\obj(60,10){$\Sigma CP^\infty_+ \vee \Sigma CP^\infty_+$}
\obj(40,0){$BGL_1K_p^\wedge$}
\mor(5,10)(15,10){}
\mor(25,10)(36,10){}
\mor(44,10)(53,10){}
\mor(40,10)(40,0){}
\mor(22,9)(36,1){$0$}
\mor(58,9)(44,1){$\widetilde{u}$}
\enddc

The top row is the cofibre sequence associated to the pushout (*). Since the map $S^1\times CP^\infty \rightarrow BGL_1K_p^\wedge$ is nullhomotopic on $CP^\infty\vee CP^\infty$, it factors through $\Sigma CP^\infty_+ \vee \Sigma CP^\infty_+$ as $\widetilde{u}$.

The map $u$ gives two units $u_1,\,u_2$ in the ${K_p^\wedge}^0(CP^\infty)$. In the Mayer Vietoris sequence for the Thom spectrum, these describe the map $K_p^\wedge \wedge CP^\infty_+ \vee K_p^\wedge \wedge CP^\infty_+ \rightarrow K_p^\wedge \wedge CP^\infty_+ \vee K_p^\wedge \wedge CP^\infty_+$ as the matrix, 
\[ \left( \begin{array}{ccc}
1 & u_2 \\
u_1 & 1  \end{array} \right)\] 

In fact, $u_1$ and $u_2$ are equal because each summand in $\Sigma CP^\infty_+$ of $\Sigma CP^\infty_+ \vee \Sigma CP^\infty_+$ is the cofibre of the map $CP^\infty \rightarrow LCP^\infty=S^1\times CP^\infty$ given by the inclusion of the constant loops and both can be defined by the same diagram,

\begindc{\commdiag}[5]
\obj(0,25){}
\obj(-25,-5){}
 \obj(0,20){$CP^\infty$}
\obj(0,10){$S^1 \times CP^\infty$}
\obj(0,0){$\Sigma CP^\infty_+$}
\obj(30,10){$BGL_1(K_p^\wedge )$}
\mor(0,20)(27,11){0}
\mor(4,10)(26,10){$\widehat{f}$}
\mor(0,20)(0,10){}
\mor(0,10)(0,0){}
\mor(3,1)(27,9){$u$}
\enddc

In terms of $u$, we can rewrite the Mayer Vietoris sequence as the long exact sequence,
$$  \ldots \rightarrow {K_p^\wedge}_*(C P^\infty) \stackrel{u-1}{\rightarrow} {K_p^\wedge}_*(C P^\infty) \rightarrow \pi_*((LC P^\infty)^{Lf})\rightarrow \ldots \ \ \ \ \ \ \  (\alpha) $$

\end{proof}

To calculate $\pi_*(THH^{K_p^\wedge}(K/p,f))$, it remains to understand the map $u$. This is done as follows:

\prop The adjoint of the map $u: \Sigma CP^\infty_+ \rightarrow BGL_1R$, is homotopy equivalent to the composite $\Sigma^2 CP^\infty_+ \stackrel{\mu}{\rightarrow}CP^\infty \stackrel{f}{\rightarrow} B^2GL_1K_p^\wedge$, where $\mu$ is the composition $\Sigma^2 CP^\infty_+ \simeq S^2 \wedge CP^\infty_+ \stackrel{\sigma \wedge id}{\rightarrow} CP^\infty \wedge CP^\infty_+ \rightarrow CP^\infty$.

\begin{proof}
The following diagram commutes:

\begindc{\commdiag}[5]
\obj(-15,15){}
\obj(0,-5){}
\obj(20,10){$S^1 \wedge LCP^\infty$}
\obj(20,0){$ CP^\infty$}
\obj(46,10){$S^1 \wedge LB^2GL_1(K_p^\wedge)$}
\obj(46,0){$B^2GL_1(K_p^\wedge )$}
\obj(0,10){$S^1 \wedge (S^1\times CP^\infty)$} 
\mor{$S^1 \wedge (S^1\times CP^\infty)$}{$S^1 \wedge LCP^\infty$}{$\simeq$}
\mor{$S^1 \wedge LCP^\infty$}{$ CP^\infty$}{$\mathit{ev}$}
\mor{$S^1 \wedge LCP^\infty$}{$S^1 \wedge LB^2GL_1(K_p^\wedge)$}{$S^1 \wedge Lf$}
\mor{$ CP^\infty$}{$B^2GL_1(K_p^\wedge )$}{$f$}
\mor{$S^1 \wedge LB^2GL_1(K_p^\wedge)$}{$B^2GL_1(K_p^\wedge )$}{$\mathit{ev}$}
\enddc 

 Consider the inclusion of the based loops $BGL_1K_p^\wedge \hookrightarrow LB^2GL_1K_p^\wedge$. Under the composite, 
$$S^1 \times BGL_1K_p^\wedge \rightarrow S^1 \times LB^2GL_1K_p^\wedge \stackrel{\mathit{ev}}{\rightarrow} B^2GL_1K_p^\wedge,$$
 the copies $S^1\times *$ and $*\times BGL_1K_p^\wedge$ map trivially. Thus, it factors through $S^1 \wedge BGL_1K_p^\wedge$ as $\Sigma BGL_1K_p^\wedge \stackrel{\sigma}{\rightarrow} B^2GL_1K_p^\wedge$. We are trying to figure out the map 
 $$S^1 \times LCP^\infty \rightarrow S^1 \times LB^2GL_1K_p^\wedge \rightarrow S^1 \times BGL_1K_p^\wedge \rightarrow B^2GL_1K_p^\wedge $$
 
 Then, this factors through 
$$S^1 \wedge LCP^\infty \rightarrow S^1 \wedge LB^2GL_1K_p^\wedge \rightarrow S^1 \wedge BGL_1K_p^\wedge \stackrel{\sigma}{\rightarrow} B^2GL_1K_p^\wedge.$$
  
Also $LCP^\infty \rightarrow BGL_1 K_p^\wedge$ factors through $S^1\wedge CP^\infty_+$ as u. Putting all the remarks together, we have a commutative diagram,

\begindc{\commdiag}[5]
\obj(-5,50){}
\obj(71,-10){}
\obj(0,40){$S^2\wedge CP^\infty_+$}
\obj(66,40){$S^1\wedge BGL_1K_p^\wedge$}
\obj(20,20){$S^1 \wedge LCP^\infty$}
\obj(20,0){$ CP^\infty$}
\obj(46,20){$S^1 \wedge LB^2GL_1(K_p^\wedge)$}
\obj(66,0){$B^2GL_1(K_p^\wedge )$}
\obj(0,20){$S^1 \wedge (S^1\times CP^\infty)$} 
\mor{$S^1 \wedge (S^1\times CP^\infty)$}{$S^1 \wedge LCP^\infty$}{$\simeq$}
\mor{$S^1 \wedge LCP^\infty$}{$ CP^\infty$}{$ev$}
\mor{$S^1 \wedge LCP^\infty$}{$S^1 \wedge LB^2GL_1(K_p^\wedge)$}{$S^1 \wedge Lf$}
\mor{$ CP^\infty$}{$B^2GL_1(K_p^\wedge )$}{$f$}
\mor{$S^1 \wedge LB^2GL_1(K_p^\wedge)$}{$B^2GL_1(K_p^\wedge )$}{$ev$}
\mor{$S^2\wedge CP^\infty_+$}{$S^1 \wedge (S^1\times CP^\infty)$}{}
\mor{$S^1\wedge BGL_1K_p^\wedge$}{$S^1 \wedge LB^2GL_1(K_p^\wedge)$}{}
\mor{$S^2\wedge CP^\infty_+$}{$S^1\wedge BGL_1K_p^\wedge$}{$\Sigma u$}
\mor{$S^1 \wedge LCP^\infty$}{$S^1\wedge BGL_1K_p^\wedge$}{$\Sigma \widehat{f}$}
\mor{$S^1 \wedge LCP^\infty$}{$S^2\wedge CP^\infty_+$}{}
\mor{$S^1\wedge BGL_1K_p^\wedge$}{$B^2GL_1(K_p^\wedge )$}{$\sigma$}
\enddc 

The left hand vertical map from $S^2\wedge CP^\infty_+$ to $S^1\wedge (S^1\times CP^\infty)$ is the inclusion of a factor in the splitting of the suspension of $S^1\wedge (S^1\times CP^\infty)\simeq (S^2\wedge CP^\infty_+)\vee(S^1\wedge CP^\infty).$

It follows that $\widetilde{u}\simeq \sigma \circ \Sigma u \simeq f\circ g$, where,
$$g : S^2\times \Sigma CP^\infty_+ \rightarrow S^1 \wedge (S^1 \times CP^\infty) \simeq S^1 \wedge LCP^\infty \stackrel{ev}{\rightarrow} CP^\infty$$
and the composition $g\simeq \mu$.
\end{proof}

\section{The Structure of $GL_1(K_p^\wedge)$}

In this section, we prove a splitting of $GL_1K_p^\wedge$ using the logarithm $l_p:gl_1 K_p^\wedge \rightarrow K_p^\wedge$ defined by Rezk (see \cite{R2}). Throughout this section, we assume that $p$ is an odd prime. 
 
\prop (Rezk, \cite{R2}) Let $R$ be an $E_\infty$ ring spectrum. Then there is a logarithmic cohomology operation, $l_{p,n}$, from $gl_1(R)$ to $L_{K(n)}(R)$ for every $n$, and prime $p$. If $R$ is $K(n)$-local, this is a map from $gl_1(R)$ to $R$. When $n=1$, $l_p:gl_1R \rightarrow R$ is given by the formula: 

$$ l_p(x) = - \frac{1}{p} log(\frac{\psi (x)}{x^p} )$$ 

[Recall that a $\theta$-algebra structure is described by operations $\psi$ and $\theta$ ($\psi$ is a ring homomorphism) such that 
$ \psi (x)= x^p + p \theta (x) $.]

\prop Suppose that $R=K_p^\wedge$. The operation $l_p:gl_1K_p^\wedge \rightarrow K_p^\wedge$ factors through $ku_p^\wedge$, the connective cover of $K_p^\wedge$. On homotopy groups, the map is an isomorphism on $\pi_n$ for $n>2$. At $n=2$, it is 0. And for $n=0$, this is the map  
$$Z_p^\times \cong Z/(p-1) \times Z_p \stackrel{p_2}{\rightarrow} Z_p$$

\begin{proof}
The spectrum $K_p^\wedge$ is $K(1)$-local, and the operation $\psi$ is the Adams operation $\psi_p$. Since $gl_1K_p^\wedge$ is connective, the map $l_p$ factors through $ku_p^\wedge$. Recall, that the homotopy groups of $gl_1K_p^\wedge$ are given by 

$$ \pi_n(gl_1K_p^\wedge) = \left\{
\begin{array}{rl}
({K_p^\wedge}^0(S^n))^\times=\pi_n(K_p^\wedge) & \text{if } n>0\\
({K_p^\wedge}^0(S^0))^\times=\pi_0(K_p^\wedge)^\times & \text{if } n=0
\end{array} \right.$$

Since $\pi_nK_p^\wedge$ is nonzero only for even n, it suffices to restrict our attention to even dimensional spheres. The $K$-theory of $S^{2n}$ is generated by $\epsilon$ where $1 - \epsilon =$ the tangent bundle of $S^{2n}$.  Hence, 
$$ \pi_{2n}(gl_1(K_p^\wedge))= {\widetilde{gl_1(K_p^\wedge)}}^0(S^{2n})= ({\widetilde{K_p^\wedge}}^0(S^{2n}))^\times= 1+ \epsilon \pi_{2n}(K_p^\wedge)$$ 
To calculate $l_p$ on $\pi_{2n}gl_1K_p^\wedge$, one needs to compute $l_p(1+k\epsilon)$ for $1+k \epsilon \in {gl_1K_p^\wedge}^0(S^{2n})=\pi_0(gl_1({K_p^\wedge}^{S^{2n}}))$. To accomplish this, we need to calculate $\psi_ p(\epsilon)$. The map $p:(S^2)^n \rightarrow S^{2n}$ which quotients out the lower cells, induces an injection in $K$-theory, and splits  $\epsilon$ as the product 
$$p^*(\epsilon)= \prod(1-L_i)$$
where $L_i$ is the canonical line bundle over the $i^{th}$ copy of $S^2=CP^1$. Since the Adams operation $\psi_p$ raises line bundles to the $p^{th}$ power,    
$$\psi_p(L_i)= L_i ^p$$ 
\begin{align}
\implies \psi_p(1-L_i) &= 1-L_i^p \nonumber \\
&= 1-(1-(1-L_i))^p \nonumber
\end{align}
 The element $1-L_i$ lies in the $K$-theory of $S^2$, so it squares to 0. Therefore,  
\begin{align}
\psi_p(1-L_i) &= 1-(1-p(1-L_i)) \nonumber \\
&= p(1-L_i) \nonumber
\end{align}
$$\implies \psi_p(\epsilon) = p^n \epsilon.$$
$$\implies \psi_p(1+\epsilon)= 1+ p^n\epsilon$$
Hence,
\begin{align}
l_p(1+k\epsilon) &= - \frac{1}{p} log(\frac{\psi (1+k\epsilon)}{(1+k\epsilon)^p} ) \nonumber \\
&= - \frac{1}{p} log(\frac{1+p^n k\epsilon}{(1+k\epsilon)^p} ) \nonumber \\
&\equiv  - \frac{1}{p}log( 1+(p^n -p) k \epsilon)\,(\mode\ p) \nonumber
\end{align}
which becomes multiplication by $1- p^{n-1}$ ($\mode\ p$) if $n > 0$. Since the homotopy group $\pi_{2n}(gl_1K_p^\wedge)=Z_p$ for $n>0$, this is an isomorphism for $n >1$. For $n=1$, this map is $0$. For $n=0$, the map $l_p : Z_p^\times \cong \mu_{p-1} \times Z_p\rightarrow Z_p$ is given by  
$$- \frac{1}{p} log(x^{1-p} )$$
This map has kernel $\nu_{p-1}$, the group of $(p-1)^{st}$ roots of unity, as it takes $p$-adic integers of the form $1+pk$ to 
\begin{align}
l_p(1+pk)&=-\frac{1}{p} log((1+pk)^{1-p}) \nonumber \\
         &= -\frac{1}{p}log(1+p(1-p)k)\nonumber \\
         &= -(1-p)k + O(p) \nonumber \\
         &\equiv -k (\mode\ p) \nonumber
\end{align}         
 Therefore, the map $l_p$ on $Z_p^\times = \nu_{p-1} \times Z_p$, has kernel $\nu_{p-1}$ and is an isomorphism onto $Z_p$ . 
\end{proof}
 
Recall that the spectrum $ku_p^\wedge$ splits into Adams summands, 
$$ku_p^\wedge \simeq B \vee \Sigma^2 B \ldots \Sigma ^{2p-4} B $$
where $B$ is the $p$-adic Adams summand($\pi_*(B)= Z_p[v_1]$). Using this, we identify the image of the logarithmic cohomology operation. We construct $K_p(\widehat{2})$ from the spectrum $ku_p^\wedge$ by killing the $2^{nd}$ homotopy group:

\defn Let $B_2$ be the 2-connective cover of $B$. Define 
$$K_p(\widehat{2})=B \vee \Sigma^2 B_2 \ldots \vee \Sigma ^{2p-4}  B $$

\prop There is a split cofibre sequence,   
$$ H\nu_{p-1} \vee \Sigma^2 HZ_p \rightarrow gl_1(K_p^\wedge) \rightarrow K_p(\widehat{2}) $$

\begin{proof}
From the definition above, note that $gl_1K_p^\wedge \stackrel{l_1}{\rightarrow}ku_p^\wedge \rightarrow K_p(\widehat{2})$ is surjective on homotopy groups. The fibre $F$ has homotopy only in dimensions 0 and 2. The Postnikov tower of $F$ then is a cofibre sequence, 
$$\Sigma^2 HZ_p \rightarrow F \rightarrow H\nu_{p-1} \rightarrow \Sigma^3 HZ_p$$

 Since the group $H^3(H\nu_{p-1};Z_p)=0$, the sequence splits and one obtains 
$$F\simeq H\nu_{p-1}\vee \Sigma^2 HZ_p$$    
Therefore, there is a cofibre sequence 
$$ H\nu_{p-1} \vee \Sigma^2 HZ_p \rightarrow gl_1(K_p^\wedge) \rightarrow K_p(\widehat{2}) $$
The next term in this sequence is  
$$\Sigma (H\nu_{p-1} \vee \Sigma^2 HZ_p)\simeq \Sigma H\nu_{p-1} \vee \Sigma^3 HZ_p$$
and the next map is $K_p(\widehat{2})\rightarrow \Sigma H\nu_{p-1} \vee \Sigma^3 HZ_p$. Since the spaces in the Adams summands are retracts of $bu_p^\wedge$, their homology concentrated in even dimensions. Therefore, 
$$[\Sigma^{2k} B, \Sigma H\nu_{p-1} \vee \Sigma^3 HZ_p] \cong H^1(B;\nu_{p-1})\oplus H^3(B;Z_p) \cong 0$$
Since the spectrum $B_2$ is 3-connected,  
$$[\Sigma^2B_2,\Sigma H\nu_{p-1} \vee \Sigma^3 HZ_p] \cong H^{-1}(B_2;\nu_{p-1})\oplus H^1(B_2;Z_p)\cong 0$$
$$\implies [K_p(\widehat{2}),H^1(B;\nu_{p-1})\oplus H^3(B;Z_p)]=0$$
Hence, the cofibre sequence splits and 
$$gl_1(K_p^\wedge) \simeq K_p(\widehat{2})\vee H\nu_{p-1} \vee \Sigma^2 HZ_p.$$ 
\end{proof}

We will use this decomposition later to calculate  homotopy classes of extensions. For that, we also have to understand how the splitting looks like when we map a space $X$ to $GL_1(K_p^\wedge)$. Recall, $[X,GL_1(K_p^\wedge)]= {K_p^\wedge}^0(X)^\times$. The map $l_p$ gives the way to map this to $[X,K_p(\widehat{2})]$. The map $ {K_p^\wedge}^0(X)^\times \rightarrow H^0(X;\nu_{p-1})$ is the composite  
$$X\rightarrow GL_1(K_p^\wedge)\rightarrow \pi_0GL_1(K_p^\wedge) \cong Z_p^\times \cong \nu_{p-1}\times Z_p \rightarrow \nu_{p-1}\simeq K(\nu_{p-1},0)$$ 

 The third factor is $\Sigma^2 HZ_p$, and we have to understand the map from $H^2(X;Z_p)$ to ${K_p^\wedge}^0(X)^\times$. Now, $H^2(X;Z_p)= [X,K(Z_p,2)]=[X,{CP^\infty}_p^\wedge]$. The space $CP^\infty$ classifies line bundles which are invertible elements in $K$-theory. 

\prop The map $H^2(X;Z_p) \rightarrow {K_p^\wedge}^0(X)^\times$ is given by $f\in [X,{CP^\infty}_p^\wedge]\mapsto L^f$ where $L^f$ is the line bundle classified by $f$.
 
 \begin{proof}
  The formula in the statement of the proposition defines a map of infinite loop spaces ${CP^\infty_p}^\wedge \rightarrow GL_1K_p^\wedge$, and hence, a map of spectra $\Sigma^2 HZ_p \rightarrow gl_1K_p^\wedge$. Composing it with $l_p$, we get 
  
  \begin{align}
  l_p(L^f)&= - \frac{1}{p} log(\frac{\psi_p (L^f)}{(L^f)^p} ) \nonumber \\
          &= -\frac{1}{p} log(\frac{(L^f)^p}{(L^f)^p}) \nonumber \\
          &= -\frac{1}{p}log(1) \nonumber \\
          &= 0 \nonumber
 \end{align}

 The computation above shows that the composition $\Sigma^2 HZ_p \rightarrow gl_1(K_p^\wedge) \rightarrow K_p(\widehat{2})$ equals 0. Therefore, it  factors through $\nu_{p-1}\times \Sigma^2 HZ_p$ in the diagram,

\begindc{\commdiag}[5]
\obj(-15,20){}
\obj(0,-5){}
\obj(20,13){$\Sigma^2 HZ_p$}
\obj(20,0){$gl_1(K_p^\wedge)$}
\obj(40,0){$ K_p(\widehat{2})$}
\obj(0,0){$H\nu_{p-1} \times \Sigma^2 HZ_p$} 
\mor(6,0)(17,0){}
\mor(23,0)(39,0){}
\mor(19,12)(1,1){}
\mor(20,13)(20,0){}
\enddc 
 
To complete this proof, we need to show that the map $\Sigma^2 HZ_p \rightarrow H\nu_{p-1}\vee \Sigma^2 HZ_p \rightarrow \Sigma^2 HZ_p$ is an equivalence. The only non zero homotopy group of $\Sigma^2 HZ_p$ is $\pi_2$, so it suffices to check that the map $[S^2,CP^\infty] \rightarrow H^2(S^2; Z_p)$ as described by the statement is an isomorphism. The left group is isomorphic to $Z_p$, via $k\mapsto L^k$, $L=$ the tangent bundle of $S^2$. The right group is $H^2(S^2; Z_p) \cong Z_p$ inside $K_p^\wedge(S^2)^\times$ as elements $1+k\epsilon$, $\epsilon = 1- L$. The map between the two is $L^k \mapsto (1-\epsilon)^k= 1-k\epsilon$ because $\epsilon^2=0$, and is evidently an isomorphism. 
 \end{proof}

\section{Calculation of $\THH$}

In this section, we complete the computation of $\THH$ for odd primes $p$. We first parameterise the homotopy classes of extensions $f$,

\begindc{\commdiag}[5]
\obj(-23,13){}
\obj(0,10){$S^2$}
\obj(20,10){$\Sigma BGL_1(K_p^\wedge)$}
\obj(0,0){$CP^\infty$}
\obj(20,0){$B^2GL_1(K_p^\wedge)$}
\mor(0,10)(15,10){$\Sigma (1-p)$}
\mor(20,10)(20,0){$\sigma$}
\mor(1,0)(15,0){$f$}
\mor(0,10)(0,0){$\sigma$}
\enddc 
using the results of the previous section.

Recall that, 
$$GL_1(K_p^\wedge)= \nu_{p-1} \times K(Z_p,2) \times \Omega^\infty K_p(\widehat{2})$$
$$\implies B^2GL_1(K_p^\wedge)= B^2 \nu_{p-1} \times K(Z_p,4)\times \Omega^\infty \Sigma^2 K_p(\widehat{2})$$
 
The condition on the map $f$ is that its restriction to $S^2$ is $1-p$. The homotopy classes of maps from $S^2$ to $B^2GL_1(K_p^\wedge)$ is split into  three factors, 

1. $[S^2,B^2\nu_{p-1}]=H^2(S^2;\nu_{p-1})\cong \nu_{p-1}$ 

2. $[S^2,K(Z_p,4)]= H^4(S^2;Z_p)=0$

3. $[S^2,\Omega^\infty \Sigma^2 K_p(\widehat{2})]= [S^2,\Omega^\infty (\Sigma^2B \vee \Sigma^4 B_2 \vee \Sigma^6 B \ldots \vee \Sigma^{2p-4}B] 
                                             = [S^2,\Omega^\infty \Sigma^2 B]
                                             =B^2(S^2)\cong Z_p $

In the splitting, 
$$[S^2,B^2GL_1(K_p^\wedge)] = \nu_{p-1} \oplus B^2(S^2)\oplus H^4(S^2; Z_p)= \nu_{p-1} \oplus Z_p \oplus 0,$$
$1-p$ is in the factor $Z_p$, where it equals $l_p(1-p)= \alpha_p$ and,
\begin{align}
\alpha_p        &= - \frac{1}{p} log((1-p)^{1-p} )\nonumber \\
                &\cong -\frac{1}{p} log ( 1- (1-p)p) \nonumber \\
                &\cong -1(\mode \ p) \nonumber
\end{align}        

\subsection{Calculation at the prime 3}

Let us begin the calculation at the prime $3$. The cofibre sequence for $gl_1K_3^\wedge$ is  
$$HZ/2 \vee \Sigma^2 HZ_3 \rightarrow gl_1(K_3^\wedge) \rightarrow  K_3(\widehat{2})$$
and 
$$K_3(\widehat{2}) = B\vee \Sigma ^2 B_2.$$
Therefore, $$GL_1K_p^\wedge = Z/2 \times K(Z_p,2)\times \Omega^\infty B \times \Omega^\infty B_2.$$

We will study the extension to $CP^\infty$ of the map $1-p$, to the four factors $Z/2,\, K(Z_3,2)$, $\Omega^\infty B,\, \Omega^\infty B_2$ one by one. Let us start with the factor $B$. The Adams summands are the eigenspaces of the action of the $(p-1)^{st}$ roots of unity by Adams operations. The spectrum $B$ is fixed by all the Adams operations. The projection from ${K_p^\wedge}^*(X)$ to $B^*(X)$ is given by 
$$\pi = \frac{1}{p-1}(1+ \psi_\zeta + \psi_{\zeta^2} +\ldots +\psi_{\zeta^{p-2}}),$$
where $\zeta \in \nu_{p-1} \subset Z_p^\times$.

For the prime $3$, we can take $\zeta = -1$ and then the projection operator is 
$$\pi =  \frac{1+\psi_{-1}}{2}$$
Let us start by working out an example. 

\ex 
Consider the element $\beta L \in {K_3^\wedge} ^2(CP^\infty)$ where $\beta$ is the Bott element. Applying the projection, we get  
\begin{align}
\pi (\beta L) &= \frac{\beta (L-L^{-1})}{2} \nonumber 
\end{align}

Restricting to $S^2$, using $L=1-\epsilon$ and $\epsilon^2= 0$, we obtain 
\begin{align}
\frac{\beta ((1-\epsilon) -(1-\epsilon)^{-1})}{2} &= \frac{\beta ((1-\epsilon)-(1+\epsilon))}{2} \nonumber \\
                                                  &= -\beta \epsilon \nonumber \\
                                                  &= -1 \nonumber
\end{align} 
                                             
In order for it to be an extension of the kind required, this restriction must be $\alpha_3$, so we multiply by $-\alpha_3$. This defines,
$$f= -\alpha_3 \frac{\beta (L-L^{-1})}{2}.$$

Recall that, $\THH^{K_3^\wedge}(K/3,f)$ is the cofibre of  
$$ K_3^\wedge \wedge CP^\infty \stackrel{u-1}{\rightarrow} K_3^\wedge \wedge CP^\infty  \ \ \ \ \ \  \ \ (\alpha)$$
where $u\in {K_3^\wedge}^0(CP^\infty) ^\times= [CP^\infty_+,GL_1(K_3^\wedge)]$ is the adjoint of,

\begindc{\commdiag}[5]
\obj(-25,-3){}
\obj(0,13){}
 \obj(0,10){$S^2 \wedge CP^\infty_+$}
\obj(0,0){$CP^\infty$}
\obj(20,0){$B^2GL_1(K_3^\wedge )$}
\mor(3,0)(16,0){$f$}
\mor(0,10)(0,0){$\mu$}
\mor(2,9)(18,1){$\widetilde{u}$}
\enddc 

 The group structure of $CP^\infty$ classifies tensor product of line bundles so, $\mu^*L=L\otimes L$. This implies, 
 \begin{align}
\mu^*(f)&=-\alpha_3 \frac{\beta(L\otimes L - L^{-1}\otimes L^{-1})}{2} \nonumber 
\end{align}                        

 The $K$ theory of $S^2$ is generated by $\epsilon = 1-L$ with $\epsilon^2 =0$. We can rewrite the equation using the generator  
\begin{align}
\mu^*(f) &=-\alpha_3 \frac{\beta((1-\epsilon)\otimes L - (1+\epsilon) \otimes L^{-1})}{2} \nonumber \\
                       &= -\alpha_3 \frac{\beta \epsilon \otimes (L + L^{-1})}{2} \nonumber
\end{align}

Using the suspension isomorphism (given by $\beta \epsilon =1$) we get, 
$$\mu^*(f)= -\alpha_3 \frac{ L + L^{-1}}{2}$$ 

To get $u$ we need to invert the logarithmic cohomology operation. Suppose that $u=h(x)\in {K_3^\wedge}^0(CP^\infty)^\times$. Then, we have to solve, 
$$- \frac{1}{3} log(\frac{\psi_3 (h(x))}{h(x)^3})= -\alpha_3 \frac{ L +L^{-1}}{2}$$
 $$\implies \frac{\psi_3 (h(x))}{h(x)^3}= exp( 3\alpha_3 \frac{ L + L^{-1}}{2})\ \ \  \ \   \ \ (*) $$

Note that $\psi_3(x)= 1- (1-x)^3$ and hence, 
$$\frac{ h(1-(1-x)^3)}{h(x)^3}= exp( 3\alpha_3 \frac{ L + L^{-1}}{2})$$

Let us look at the equation $(\mode \ 3^2,x^3)$. The right side of the equation can be written in terms of $x$ using $L= 1-x$, and then, $L^{-1}= 1+x+x^2 \, (\mode \ 3^2,\, x^3)$. Therefore, the right side simplifies to  
\begin{align}
exp(3\alpha_3 \frac{ L + L^{-1}}{2})&= exp( 3 \alpha_3 \frac{2 + x^2}{2}) \nonumber \\ 
                                    &=1 + 3\alpha_3 + 3\frac{\alpha_3 x^2}{2} \nonumber
\end{align}

Now we will simplify the left side of (*). Suppose that $h(x)= a + bx+ cx^2$. In order to solve the equation, we have to invert $l_3$. We know that $l_3$ has a kernel $Z/2\vee K(Z_3,2)$, so the equation can be solved once we know the restriction to these.

In the part of $HZ/2$,  $\sigma:S^2 \rightarrow CP^\infty$ induces an isomorphism in $H^2(-; Z/2)$. Therefore, the extension is $0$ here. The map ${K_3^\wedge}^0(CP^\infty)^\times \rightarrow H^0(CP^\infty;Z/2)$ sends $a \mapsto a\,(\mode\ 3)$ (identifying $Z/2$ with the group of units in $\mathbb{F}_3$). Therefore, since $\mu^*(0)=0$, we get the equation   
$$a\equiv 1\,(\mode\ 3)$$

In the factor $K(Z_3,2)$, there is no restriction on $f$. Assume that it is trivial, so $\mu^*(0)=0$.  This maps into $GL_1(K_3^\wedge)$ by taking a line bundle over $CP^\infty$ to the corresponding unit in $K$-theory. If we look at $k\in Z_p= H^2(CP^\infty; Z_3)=[CP^\infty,K(Z_3,2)]$, this is the line bundle $L^k= (1-x)^k= 1-kx + \frac{k(k-1)}{2}x^2 \,(\mode\, x^3)$. This is the only factor that gives a non zero coefficient of $x$ so, we get that $b=0$.

Therefore, $h(x)= a + c x^2(mode\ 3^2,\,x^3)$ and $a\equiv 1(\mode \ 3)$. The left side of (*) is 
\begin{align}
\frac{\psi (h(1-(1-x)^3))}{h(x)^3}&\equiv \frac{h(3x-3x^2+x^3)}{h(x)^3} \nonumber \\
                                  &\equiv \frac{a}{a^3 + 3ca^2x^2} \nonumber \\
                                  &\equiv a^{-2}( 1- 3\frac{c}{a}x^2)\,(\mode\ 3^2,\,x^3) \nonumber 
\end{align}

Working $(\mode\ 3^2,\, x^3)$, we have  
$$a^{-2}( 1- 3\frac{c}{a}x^2)= 1 + 3\alpha_3 + 3\frac{\alpha_3 x^2}{2}$$
$$\implies a\equiv 1+3\alpha_3 (\mode\ 3^2) \ \ and \ \ c\equiv \alpha_3 (\mode\ 3)$$
Thus $a= 1+ 3.(unit)$ and $c$ is a unit (since $\alpha_3$ is a unit). 

Therefore, $u-1$ looks like $3.(\mathit{unit})+ x^2(\mathit{unit})$. We can choose a different parameterisation for $K$-theory  of $CP^\infty$ to assume that $u-1= 3+ x^2$.

Now ${K_3^\wedge}_*(CP^\infty)= {K_3^\wedge}_*\left\{\beta_0,\beta_1,\ldots \right\}$ where $\beta_i$ is dual to $x^i$. Therefore, 
$$ < (u-1)(\beta_i),x^j>= <\beta_i, x^j(3+ x^2)>= \left\{
\begin{array}{rl}
3 & \text{if } j=i\\
1 & \text{if } j=i-2\\
0 & \text{} otherwise
\end{array} \right.$$
Therefore, the map $u-1$ on ${K_p^\wedge}_*(CP^\infty)$ is given by 
$$(u-1)(\beta_i)=\left\{
 \begin{array}{rl}
3\beta_i & \text{if } i=0,1\\
3\beta_i + \beta_{i-2} & \text{if } i>1
\end{array}\right.$$
 
Following the cofibre $(\alpha)$, we understand that $u-1$ is injective, and its cokernel has two copies of $Z/(3^\infty)$ in even dimensions. Thus, 
$$\pi_k(\THH^{K_3^\wedge}(K/3),f)=\left\{
 \begin{array}{rl}
0  & \text{if $k$ is odd }\\
Z/(3^\infty) \oplus Z/(3^\infty) & \text{if $k$ is even}
\end{array}\right.$$
completing the calculation in this example.

\begindc{\commdiag}[5]
\obj(0,0){}
\obj(0,5){}
\enddc

 Now we perform the calculation at the prime 3 for all extensions that are non trivial only on the factor $\Omega^\infty B$ of $GL_1(K_3^\wedge)$. The extension in the example was of this kind. So, we are looking at elements in $B^2(CP^\infty)$ which restrict to $\alpha_3$ in $S^2$. 

An element in ${K_3^\wedge}^2(CP^\infty)$ is given by $\beta g(x)$. Therefore, an element in $B^2(CP^\infty)$ is 
$$\pi(\beta g(x))= \frac{\beta(g(x)-g(1-\frac{1}{1-x}))}{2}$$

Suppose that $g(x)= a'+b'x+c'x^2 \,(\mode\ 3^2,\, x^3)$. Restricting to $S^2$ (using $x=\epsilon$ and $\epsilon^2= 0$) we get $b'$. We need to get $\alpha_3$. Thus, to get an extension we must have  $b'= \alpha_3$. This gives us all possible extensions $f$ on the factor $B$. Let us work as before $(\mode \ 3^2,\,x^3)$. Then, 

\begin{align}
f &= \frac{\beta(g(x)-g(1-\frac{1}{1-x}))}{2} \nonumber \\
                &= \frac{\beta(a'+ b'x+c'x^2 - g(-x -x^2))}{2} \nonumber \\
                &= \frac{\beta(a'+ b'x+ c'x^2-(a'-b'x-b'x^2+c'x^2))}{2}\nonumber \\
                &=\frac{\beta(2b'x +b'x^2)}{2} \nonumber \\
                &=\beta b'x + \frac{\beta b'}{2} x^2\nonumber
\end{align}                

We have to calculate $u$ using 

\begindc{\commdiag}[5]
\obj(-25,-3){}
\obj(0,13){}
 \obj(0,10){$S^2 \times CP^\infty$}
\obj(0,0){$CP^\infty$}
\obj(20,0){$B^2GL_1(K_3^\wedge )$}
\mor(3,0)(16,0){$f$}
\mor(0,10)(0,0){$\mu$}
\mor(2,9)(18,1){$u$}
\enddc

By definition, the multiplication map takes $x$ to the formal group, which for $K$-theory is the multiplicative group. Therefore,
$$x\mapsto \epsilon \otimes 1 + 1\otimes x -\epsilon \otimes x$$  
\begin{align}
\implies x^2 &\mapsto (\epsilon \otimes 1 + 1\otimes x -\epsilon \otimes x)^2 \nonumber \\
     &= 1\otimes x^2 + 2 \epsilon \otimes x -2\epsilon \otimes x^2 \nonumber 
\end{align}
To get $\mu^*$ we must project onto the factor $S^2\wedge CP^\infty_+$. Thus, we obtain  
$$\mu^*(x)= \epsilon \otimes 1 - \epsilon \otimes x  \ , \   \mu^*(x^2)= 2\epsilon \otimes x - 2\epsilon \otimes x^2.$$
 
Using these formulae and the suspension isomorphism $\beta \epsilon = 1$ we calculate $\mu^*(f)$.

\begin{align}
\mu^*(f) &= \beta b'(\epsilon \otimes 1- \epsilon \otimes x)+ \frac{\beta b'}{2}(2 \epsilon \otimes x-2\epsilon \otimes x^2) \nonumber \\
                       &= b'(1-x) + \frac{b'}{2}(2x-2x^2) \nonumber \\
                       &= b' - b'x^2 \nonumber
\end{align}

To get $u$, we have to invert the logarithmic cohomology operation $l_3$, as in the example. Suppose that $u= h(x)$. Then, we need to solve 
$$l_3(u)= \frac{\psi_3 (h(x))}{h(x)^3}= exp(-3b'(1-x^2))$$

We have the formula $\psi_3(x)= 1-(1-x)^3$. Similar to the example, we assume that in our extension the contribution from $HZ/2$ is $1$ and $HZ_3$ is $0$. In the same way, this implies that if $h(x)= a + bx +c x^2$ ,
$$a\equiv 1(\mode \ 3) \ , \ b=0$$  
Then, the equation becomes, 
\begin{align}
a^{-2}( 1- 3\frac{c}{a}x^2) &= exp(-3b'(1-x^2)) \nonumber \\
                            &= 1 - 3b' + 3b'x^2 \nonumber 
\end{align}

In the same way, we understand that the unit $u= 1+ 3.\mathit{unit} + x^2.\mathit{unit}$, and so, we obtain the same computation  

$$\pi_k(\THH^{K_3^\wedge}(K/3),f)=\left\{
 \begin{array}{rl}
0  & \text{if $k$ is odd }\\
Z/(3^\infty) \oplus Z/(3^\infty) & \text{if $k$ is even}
\end{array}\right.$$

 Now we want to see what happens if we allow extensions with non trivial contributions from the other 3 factors of $GL_1(K_3^\wedge)= Z/2 \times K(Z_3,2) \times \Omega^\infty B \times \Omega^\infty B_2$. In the part $Z/2$, the restriction $H^2(CP^\infty;Z/2)\rightarrow H^2 (S^2;Z/2)$ is an isomorphism. So, this factor always contributes trivially. 
 
 For the factor $K(Z_3,2)$, the group $[S^2,B^2K(Z_3,2)]=[S^2,K(Z_3,4)]=H^4(S^2;Z_3) =0$. Therefore, there is no condition on $f$ here. The group $H^4(CP^\infty;Z_3)$ is generated by $x^2$ and $f$ is given by $ax^2$ for some $a \in Z_p$. To compute $u$, consider:

\begindc{\commdiag}[5]
\obj(-25,-3){}
\obj(0,13){}
 \obj(0,10){$S^2 \times CP^\infty$}
\obj(0,0){$CP^\infty$}
\obj(20,0){$K(Z_3,4)$}
\mor(3,0)(16,0){$f$}
\mor(0,10)(0,0){$\mu$}
\mor(2,9)(18,1){$u$}
\enddc

Note that in this case, $\mu^*(x)= \epsilon \otimes 1 + 1\otimes x$, which implies
\begin{align}
 \mu^*(x^2) &=(\epsilon \otimes 1 + 1\otimes x)^2 \nonumber \\
                    &= 2\epsilon \otimes x + 1 \otimes x^2 \nonumber 
\end{align}                    
                     
To get $u$ we have to project to $S^2\wedge CP^\infty_+$ and apply the suspension isomorphism. Then, we get $2a x \in H^2(CP^\infty; Z_3)$. 
Recall from the previous section that, from this we get the unit by taking $L^{2a}$, where $L= (1-x)$ is the canonical line bundle. Therefore, the contribution to $u$ from this factor is $(1-x)^{2a}$.

Now if $a$ is divisible by $3$ then, we still get that our $u= 1 + 3.\mathit{unit} + x^2.\mathit{unit}$ which results in the same calculation for $\pi_*(\THH^{K_3^\wedge}( K/3,f))$. If $a$ is not divisible by $3$ then, it is a unit, so that $u= 1+ 3.\unit + x.\unit$. Therefore, by reparameterising we can write $u-1= 3 + x$. 
 $$\implies < (u-1)(\beta_i),x^j>= <\beta_i, x^j(3 +x )>= \left\{
\begin{array}{rl}
3 & \text{if } j=i\\
1 & \text{if } j=i-1\\
0 & \text{} otherwise
\end{array} \right.$$
$$\implies (u-1)(\beta_i)=\left\{
 \begin{array}{rl}
3\beta_i & \text{if } i=0\\
3\beta_i + \beta_{i-1} & \text{if } i>0
\end{array}\right.$$
 
Therefore, in this case,

$$\pi_k(THH^{K_3^\wedge}(K/3),f)=\left\{
 \begin{array}{rl}
0  & \text{if $k$ is odd }\\
Z/(3^\infty) & \text{if $k$ is even}
\end{array}\right.$$

Now consider the factor $\Sigma ^2 B_2$. We know that this is $5$-connected. So, if we look at extensions we know that they always restrict to $0\in K_3^\wedge (CP^2)$. Since we are working $(\mode\ x^3)$, this means that these extensions always give $0$.

Therefore, we get that, depending on $f$ either $ \pi_*(\THH^{K_3^\wedge}( K/3,f))= (Z/(3^\infty))^i$ in even degrees where $i=1$ or $2$ depending on $f$.

This finishes our calculation at the prime 3. 

\subsection{Calculation at primes $\geq 5$} 
Let us now look at the other odd primes and work $(\mode \ x^p,\,p^2)$. Recall that there is a splitting  
$$GL_1(K_p^\wedge)= \nu_{p-1} \times K(Z_p,2)\times \Omega^\infty K_p(\widehat{2})$$ 
$$K_p(\widehat{2})=B \vee \Sigma^2 B_2\vee \ldots \Sigma ^{2p-4} B $$

We start by working in the factor $B$ of $K_p(\widehat{2})$. The projection operator from ${K_p^\wedge} ^*(X)$ to $B^*(X)$ is given by 
$$ \pi = \frac{1+ \psi_\zeta + \psi_{\zeta^2} + \ldots \psi_{\zeta^{p-2}}}{p-1}$$

Define $\kappa$ to be the composite, 
$${K_p^\wedge} ^*(CP^\infty) \stackrel{\mu^*}{\rightarrow} {K_p^\wedge} ^*(S^2 \wedge CP_+^\infty) \stackrel{\cong}{\rightarrow} {K_p^\wedge} ^{*-2} (CP^\infty)$$

First observe that the following diagram commutes:

\begindc{\commdiag}[5]
\obj(-20,15){}
\obj(0,-5){}
\obj(20,10){${K_p^\wedge}^2(CP^\infty)$}
\obj(0,10){${K_p^\wedge}^2(CP^\infty)$}
\obj(0,0){${K_p^\wedge}^0(CP^\infty)$}
\obj(20,0){${K_p^\wedge}^0(CP^\infty)$}
\mor(5,10)(15,10){$\psi_a$}
\mor(0,10)(0,0){$\kappa$}
\mor(5,0)(15,0){$\psi_a$}
\mor(20,10)(20,0){$\kappa$}
\enddc 

This implies all Adams operations hence $\pi$, commutes with $\kappa$. 

Write $x=1-L$ for the generator in ${K_p^\wedge}^2(CP^\infty)$ amd $\epsilon$ its restriction to $S^2$. We have to look for $f$ as in the diagram, 

\begindc{\commdiag}[5]
\obj(-23,13){}
\obj(0,10){$S^2$}
\obj(20,10){$\Sigma BGL_1(K_p^\wedge)$}
\obj(0,0){$CP^\infty$}
\obj(20,0){$B^2GL_1(K_p^\wedge)$}
\mor(0,10)(15,10){$\Sigma (1-p)$}
\mor(20,10)(20,0){$\sigma$}
\mor(1,0)(15,0){$f$}
\mor(0,10)(0,0){$\sigma$}
\enddc 

Suppose that  $f$ is given by $\pi (\beta g(x))$ where $g(x)= a_0 +a_1 x+\ldots +a_{p-1}x^{p-1}\,(\mode\ x^p,\,p^2)$. 

\textbf{Claim:}
$$\kappa(\beta g(x))= g'(x)(1-x)$$ 

\begin{proof}
It is enough to check this on the generators $x^n$. The multiplication takes $x$ to the formal group of $K$-theory, which is the multiplicative formal group. 
$$\mu ^*(x)= \epsilon \otimes 1 - \epsilon \otimes x$$
Therefore, 
\begin{align}
\mu ^*(\beta x^n)&= \beta(\epsilon \otimes 1 + 1\otimes x - \epsilon \otimes x)^n \nonumber \\
                 &= \beta(1 \otimes x^n + n \epsilon \otimes x^{n-1} - n \epsilon \otimes x^n) \nonumber
\end{align}                 

$\kappa$ is obtained by projecting this onto the factor $S^2\wedge CP^\infty_+$ of the product, and then applying the suspension isomorphism ($\beta\epsilon =1$). Therefore, we obtain 

\begin{align}
\kappa(\beta x^n)&= nx^{n-1} - n x^n \nonumber \\
                 &= nx^{n-1}(1-x) \nonumber \\
                 &= (x^n)'(1-x) \nonumber
\end{align}

\end{proof}

If we restrict $f$ to $S^2$, we get 
\begin{align}
\pi (\beta g(\epsilon)) &= \pi (\beta( a_0 +a_1 \epsilon))\nonumber \\
                        &= (\frac{1+ \psi_\zeta + \psi_{\zeta^2} + \ldots \psi_{\zeta^{p-2}}}{p-1})(\beta(a_0 + a_1 \epsilon ))\nonumber
\end{align}                        
The action of the Adams operations on the Bott element and $\epsilon$ are given by 
$$\psi_a(\beta)=\frac{\beta}{a} \ , \ \psi_a(\epsilon)=1-(1-\epsilon)^a= a\epsilon$$
Therefore,
\begin{align}
\pi (\beta g(\epsilon)) &= (\frac{1+ \psi_\zeta + \psi_{\zeta^2} + \ldots \psi_{\zeta^{p-2}}}{p-1})(\beta(a_0 + a_1 \epsilon ))\nonumber \\
                        &= \frac{\beta( a_0(1+ \zeta^{-1} +\zeta^{-2}+\ldots +\zeta^{2-p}) +(p-1)a_1 \epsilon) }{p-1} \nonumber \\
                        &= \frac{\beta( a_0(1+ \zeta^{-1} +\zeta^{-2}+\ldots +\zeta^{2-p}))}{p-1} +a_1 \epsilon \nonumber 
\end{align}                        
Since $\zeta$ is a $(p-1)^{st}$ root of unity, we get
\begin{align}
1+ \zeta^{-1} +\zeta^{-2}+\ldots +\zeta^{2-p} &= \zeta^{p-1} + \zeta^{p-2} +\ldots +\zeta \nonumber \\
                                              &= 0 \nonumber
\end{align}                                              
This shows that $\pi(\beta g(x))$ restricts to $a_1 \in B^2(S^2)$. Thus, we have that $a_1= l_p(1-p)= \alpha_p$.

We need to calculate $u$ from the extension $\pi(\beta g(x))$ by solving,

\begin{align}
l_p(u) &= \kappa \pi(\beta g(x))\nonumber \\
       &= \pi \kappa(\beta g(x))\nonumber \\
       &= \pi(g'(x)(1-x)) \nonumber
\end{align}

Suppose that $h(x)= g'(x)(1-x)= c_0+c_1 x +\ldots + c_{p-1}x^{p-1} (\mode \ x^p,\,p^2)$. Then,

\begin{align}
\pi (h(x)) &= \frac{1+ \psi_\zeta + \psi_{\zeta^2} + \ldots \psi_{\zeta^{p-2}}}{p-1}(h(x))\nonumber \\
           &= \sum_{i=0}^{p-2} \frac{ h(1-(1-x)^{\zeta^i})}{p-1} \nonumber
\end{align}

Let us look at the coefficient of $x^a$ in the above equation. 

\begin{align}
[\pi (x^n)]_a &= [\sum_{i=0}^{p-2} \frac{ (1-(1-x)^{\zeta^i})^n}{p-1}]_a \nonumber \\
              &=\sum_{i=0}^{p-2} \frac{[(1-(1-x)^{\zeta^i})^n]_a}{p-1} \nonumber \\
              &=\sum_{i=0}^{p-2} \frac{[ \sum {{n \choose l} (-1)^l(1-x)^{l\zeta ^i}}]_a}{p-1} \nonumber \\
              &=\sum_{i=0}^{p-2} \frac{\sum { {n\choose l} (-1^l) {{l\zeta^i}\choose a}}}{p-1} \nonumber
\end{align}               

 Since $\zeta$ is a $(p-1) ^{st}$ root of unity,
$$(\zeta) ^i + (\zeta^2)^i +\ldots + (\zeta^{p-1})^i =  \left\{
 \begin{array}{rl}
0 & \text{if } i=1,2,\ldots ,p-2\\
p-1  & \text{if } i=0,p-1
\end{array}\right.$$
The binomial coefficient ${y \choose a}$ is a polynomial in $y$ of degree $a$ with the constant term 0  and the top coefficient $1/a!$. Therefore,  
$$  {{l(\zeta)}\choose a} + {{l(\zeta^2)}\choose a} +\ldots + {{l(\zeta^{p-1})}\choose a} =  \left\{
 \begin{array}{rl}
0 & \text{if } a=1,2,\ldots ,p-2\\
\frac{p-1}{(p-1)!}l^{p-1}  & \text{if } a=p-1\\
p-1 &\text{if } a=0
\end{array}\right.$$ 
Therefore, we get 
$$[\pi (x^n)]_a = \left\{ 
\begin{array}{rl}
0 & \text{if } a=1,2,\ldots ,p-2\\
\frac{1}{(p-1)!}\sum{{n\choose l}(-1)^l l^{p-1}}  & \text{if } a=p-1\\
\sum{{n\choose l}(-1)^l} & \text{if } a=0
\end{array}\right.$$
$$\implies [\pi (x^n)]_0 = \sum{{n\choose l}(-1)^l} =\left\{
 \begin{array}{rl}
(1-1)^n = 0  & \text{if } n > 0 \\
1  & \text{if } n=0
\end{array}\right.$$
The other possible non zero coefficient is $[\pi(x^n)]_{p-1}$. If $n=0$, this must be 0. If $n>0$ this gives 
 \begin{align}
 [\pi (x^n)]_{p-1} &\equiv \frac{1}{(p-1)!}\sum{{n\choose l}(-1)^l l^{p-1}}\nonumber \\
                   &\equiv \frac{1}{(p-1)!}\sum{{n\choose l}(-1)^l l^{p-1}}\nonumber \\
                   &\equiv - \sum{{n\choose l}(-1)^l} \nonumber \\
                   &\equiv -(1-1)^n +1 \nonumber \\
                   &\equiv 1\, (\mode\ p) \nonumber
 \end{align}
 Summarising the calculation $(\mode\ p)$, we get 
$$[\pi (x^n)]_a =  \left\{ 
\begin{array}{rl}
1 & \text{if } a=0,\,n=0\\
1  & \text{if } a=p-1,\,n>0\\
0 & \text{otherwise }  
\end{array}\right.$$
 Now we are in a position to calculate $\pi(h(x))\,(\mode\ p)$ 
 \begin{align}
 \pi(h(x)) &=\pi(c_0 + c_1 x +\ldots +c_{p-1}x^{p-1}) \nonumber \\
           &=c_0\pi(1)+c_1\pi(x)+\ldots +c_{p-1}\pi(x^{p-1}) \nonumber \\
           &=c_0+ c_1 x^{p-1} + \ldots +c_{p-1}x^{p-1} \nonumber \\
           &=c_0+bx^{p-1} \nonumber
 \end{align}            
 where $c_0=a_1$ and 
 \begin{align}
 b&= c_1+ \ldots c_{p-1} \nonumber \\
  &=a_1 -2a_2 + 2a_2 -3a_3 \ldots -(p-1)a_{p-1}+(p-1) a_{p-1} -p a_p\nonumber \\
  &\equiv a_1 \,( \mode \ p) \nonumber 
 \end{align} 
Thus the equation for $u \, (\mode\ p)$ reduces to  
$$l_p (u)= a_1+ b x^{p-1}\ (\mode\ p)$$
$$\implies -\frac{1}{p}log(\frac{\psi_p(u(x))}{u^p})= a_1 + bx^{p-1}\ (\mode\ p)$$
$$\implies \frac{\psi_p(u(x))}{u^p}= exp( -p(a_1 + bx^{p-1}))= 1 -p a_1 + pb x^{p-1}(\mode\ p^2)$$
We are looking at extensions which are non trivial only on the factor $B$. This implies $u(x) \in B^0(CP^\infty)$ which implies $u$ is in the image of $\pi$. By the calculations above, this implies that $u(x)= d_0 + d_1 x^{p-1} \,(\mode\, x^p)$. Then,
 
\begin{align}
\frac{\psi_p(u(x))}{u^p}&= \frac{ d_0}{ {d_0}^p + p d_0^{p-1} d_1 x^{p-1}} \nonumber \\
                       &= (d_0)^{1-p}(1- p \frac{d_1}{d_0} x^{p-1})\nonumber
\end{align}                       

Therefore, we obtain 

$$d_0 ^{1-p}= 1- pa_1$$ 
 \begin{align}
 \implies d_0 &= (1-pa_1)^{\frac{1}{1-p}} \nonumber \\
              &= 1- \frac{p}{1-p} a_1 \nonumber \\
              &= 1- p a_1 \,(\, \mode\, p^2) \nonumber
 \end{align}             
\begin{align}
\implies d_1 &= - {d_0}^p b \nonumber \\
             &\equiv -1\, (\mode\ p) \nonumber
\end{align}              

Therefore, $d_0 = 1+ p.\unit$ and $d_1=\unit$. Thus, $u= 1+ p.\unit + \unit. x^{p-1}$. We can reparameterise so that $u= 1+ p + x^{p-1}$. 
\begin{align}
< (u-1)(\beta_i),x^j>&= <\beta_i, x^j(p +x^{p-1} )> \nonumber \\
                     &= \left\{
\begin{array}{rl}
p & \text{if } j=i\\
1 & \text{if } j=i-(p-1)\\
0 & \text{} otherwise
\end{array} \right. \nonumber
\end{align}
$$\implies (u-1)(\beta_i)=\left\{
 \begin{array}{rl}
p\beta_i & \text{if } i=0\\
p\beta_i + \beta_{i-(p-1)} & \text{if } i>0
\end{array}\right.$$

Inputting this in the long exact sequence $(\alpha)$, we get
$$ \pi_k(\THH^{K_p^\wedge}(K/p,f))=\left\{
 \begin{array}{rl}
0  & \text{if $k$ is odd }\\
(Z/(p^\infty))^{p-1} & \text{if $k$ is even}
\end{array}\right.$$
   Now let us look at what happens if we allow non trivial extensions in the other factors. Under restriction to $S^2$, $H\nu_{p-1}^2(S^2) \cong H\nu_{p-1}^2(CP^\infty) = \nu_{p-1}$. The element $1-p$ gives $1 \in \nu_{p-1}$. So, this part always contributes trivially. 

The factor $\Sigma^2 B_2$ is $(2p-1)$-connected. So, $[ CP^{p-1} , \Sigma ^2 B_2 ] = 0$. Thus $(\mode \  x^p)$ this factor is always trivial. 

Next lets look at the factor $\Sigma^2 HZ_p$. Since $[S^2,\Sigma^4 HZ_p] = HZ_p^4(S^2)= 0$, we have no condition on $f$ from this factor. The group $HZ_p ^4(CP^\infty)$ is generated by $x^2$. Suppose that $f$ is given by $ax^2 \in HZ_p ^4(CP^\infty)$. To compute the contribution to  $u$, we have the diagram,

\begindc{\commdiag}[5]
\obj(-25,-3){}
\obj(0,13){}
 \obj(0,10){$S^2 \times CP^\infty$}
\obj(0,0){$CP^\infty$}
\obj(20,0){$K(Z_p,4)$}
\mor(1,0)(18,0){$f$}
\mor(0,10)(0,0){$\mu$}
\mor(2,9)(18,1){$u$}
\enddc 

Under $\mu$, $x$ pulls back to the formal group and thus,
\begin{align}
\mu^*(x^2)&=(\epsilon \otimes 1 + 1\otimes x)^2 \nonumber \\
          &= 2\epsilon \otimes x + 1 \otimes x^2\nonumber
\end{align}           
Projecting this to the factor $S^2\wedge CP^\infty_+$, and applying the suspension isomorphism we get $2a x \in HZ_p^2(CP^\infty)$. The map from $HZ_p^2(CP^\infty) \rightarrow [CP^\infty, GL_1(K_p^\wedge)] = {K_p^\wedge} ^0(CP^\infty) ^\times $ is given by $\alpha x \rightarrow (1+x)^\alpha $. 

Therefore, if $a$ is divisible by $p$ then we still get that $u= 1 + p.\unit + x^{p-1}.\unit$. This does not change the calculation of $\THH^{K_p^\wedge}( K/p,f)$. If $a$ is not divisible by $p$, then it is a unit. Then, $u= 1+ p.\unit + x.\unit$. This can be reparameterised to $u= 1+p+x$. Then 

\begin{align}
< (u-1)(\beta_i),x^j>&= <\beta_i, x^j(p+ x)> \nonumber \\
                     &= \left\{
\begin{array}{rl}
p & \text{if } j=i\\
1 & \text{if } j=i-1\\
0 & \text{} otherwise
\end{array} \right. \nonumber 
\end{align}

$$\implies (u-1)(\beta_i)=\left\{
 \begin{array}{rl}
p\beta_i & \text{if } i=0\\
p\beta_i + \beta_{i-1} & \text{if } i>0
\end{array}\right.$$
 
Therefore, we obtain  

$$ \pi_k(\THH^{K_p^\wedge}(K/p,f))=\left\{
 \begin{array}{rl}
0  & \text{if $k$ is odd }\\
Z/(p^\infty) & \text{if $k$ is even}
\end{array}\right.$$

The other factors are $\Sigma^{2k} B$ for $k=2,3,\ldots p-2$. These correspond to the eigenspaces of the action of the Adams operations where $\psi_{\zeta^i}$ acts as $\zeta^{ki}$. The projection operator is given by 
$$\pi_k= \frac{1+\zeta^{-k} \psi_\zeta + \zeta^{-2k} \psi_{\zeta^2} +\ldots + \zeta^{-k(p-2)} \psi_{\zeta^{p-2}}}{p-1}$$

The group $[S^2, \Omega^\infty \Sigma^2\Sigma^{2k} B]= B^{2k+2}(S^2)=0$, so, there is no condition on restriction to $S^2$. Then, we may choose any $\pi_k(\beta h(x))$ for $f$, and $u$ must satisfy, 
\begin{align}
l_p(u)&= \kappa (\pi_k(\beta h(x)) ) \nonumber \\
      &=\pi_k(\kappa (\beta h(x))) \nonumber \\
      &=\pi_k(f'(x)(1-x)) \nonumber
\end{align}      

Now assume $g(x)= h'(x)(1-x)= c_0+c_1 x +\ldots + c_{p-1}x^{p-1}$. Then 
\begin{align}
\pi_k (g(x))&= \frac{1+ \zeta^{-k}\psi_\zeta + \zeta^{-2k}\psi_{\zeta^2} + \ldots +\zeta^{-k(p-2)}\psi_{\zeta^{p-2}}}{p-1}(g(x)) \nonumber \\
            &= \frac{1+ \zeta^{-k}g(1-(1-x)^\zeta) + \zeta^{-2k}g(1-(1-x)^{\zeta^2}) + \ldots \zeta^{-k(p-2)}g(1-(1-x)^{\zeta^{p-2}})}{p-1}\nonumber
\end{align}            

The following proposition is useful to complete the calculation 
\prop 
There is a polynomial $f_k(x)= x^k + a_{k+1} x^{k+1} + \ldots $ such that, $Im(\pi_k)$ has polynomials that are multiples of $f_k\,(\mode\ x^p)$ .    

\begin{proof}
These polynomials are in the $p$-adic $K$-theory of $CP^\infty$. By looking $(\mode\ x^p)$, we are restricting to the $K$-theory of $CP^{p-1}$. It splits into eigenspaces 
$${K_p^\wedge}^0(CP^{p-1})= \oplus_{k=0}^{p-2} \Lambda_k$$
  where $\Lambda_k=[CP^{p-1}, \Omega^\infty \Sigma^{2k} B]=B^{2k}(CP^{p-1})$ is the eigenspace on which the Adams operations $\psi_\zeta$ act as multiplication by $\zeta^k$. $\pi_k$ is the projection on to the eigenspace $\Lambda_k$. In this decomposition, $dim(\Lambda_0)=2$ and $dim(\Lambda_k)=1$ for all $k\geq 1$. Therefore, $\Lambda_k = span(f_k)$ for some polynomial $f_k$. To see how the polynomial $f_k$ looks like we compute $\pi_k(x)$ (note $\pi_k(1)=0$, so we don't get any information out of it). 
\begin{align}
\pi_k (x)&= \frac{1+ \zeta^{-k}\psi_\zeta + \zeta^{-2k}\psi_{\zeta^2} + \ldots +\zeta^{-k(p-2)}\psi_{\zeta^{p-2}}}{p-1}(x) \nonumber \\
            &= \frac{1}{p-1}\sum_{i=0}^{p-2}\zeta^{-ik}(1-(1-x)^{\zeta^i}) \nonumber \\
            &= \frac{1}{p-1}\sum_{i=0}^{p-2}\zeta^{-ik} \sum_{n=1}^{\infty} (-1)^{i-1}{\zeta^i \choose n} x^n \nonumber \\
            &= \frac{1}{p-1}\sum_{i=0}^{p-2}\sum_{n=1}^{\infty}   (-1)^{i-1}\zeta^{-ik}{\zeta^i \choose n} x^n \nonumber  
\end{align}    

Lets look at the coefficient of $x^n$ in the above formula. ${y \choose n}$ is a polynomial of degree $n$ in $y$, and therefore, $-1)^{i-1}y^{-k}{y \choose n}$ has terms of degree $-k$ to $-k+n$. So, if we sum the series, it is 0 if $k< n$. Thus, the first possible non zero coefficient of $x$ is in degree $k$. The coefficient of $x^k$ in $\pi_k(x)$ is given by 
\begin{align}
[\pi_k(x)]_k &\equiv \frac{1}{p-1}\sum_{i=0}^{p-2} (-1)^{i-1}\zeta^{-ik}{\zeta^i \choose k}  \nonumber \\
             &\equiv \frac{1}{p-1}\sum_{i=0}^{p-2} (-1)^{i-1}\frac{1}{k!} \nonumber \\ 
             &\equiv \frac{1}{(p-1) k!} \nonumber \\
             &\neq 0 \,(\mode\ p) \nonumber
\end{align}
               
So, this is a unit in $Z_p$. Therefore, $Im(\pi_k)= Span(f_k)$ where $f_k$ looks like $x^k+ O(x^{k+1})$.

\end{proof}

Therefore, $\pi_k(g(x))= cf_k(x)$ for some constant $c$. The equation for $u$ is 
$$l_p(u)= \pi_k(g(x))= c f_k(x)$$
$$\implies \frac{\psi_p(u)}{u^p} = exp(-p c f_k(x) )$$

If $c$ is divisible by $p$, then $(\mode\ p^2)$ the above equation is 0. If $c$ is not divisible by $p$, then the coefficient of $x^k$ in the right side is $p$ times an unit. We can solve for $u$ as in the cases before. From here, we get a contribution = $\unit.x^k$. Therefore, the unit becomes $u= 1+p.\unit + x^k.\unit$. As before, we have the long exact sequence $(\alpha)$
$$  {K_p^\wedge}_*(C P^\infty) \stackrel{u-1}{\rightarrow} {K_p^\wedge}_*(C P^\infty) \rightarrow \pi_*(\THH^{K_p^\wedge}(K/p,f)) $$ 
and,
\begin{align}
< (u-1)(\beta_i),x^j>&= <\beta_i, x^j(p+ x^k)> \nonumber \\
                     &= \left\{
\begin{array}{rl}
p & \text{if } j=i\\
1 & \text{if } j=i-k\\
0 & \text{} otherwise
\end{array} \right. \nonumber 
\end{align}
$$\implies (u-1)(\beta_i)=\left\{
 \begin{array}{rl}
p\beta_i & \text{if } i=0\\
p\beta_i + \beta_{i-k} & \text{if } i>0
\end{array}\right.$$
 
Therefore, we obtain that 
$$\implies \pi_n(\THH^{K_p^\wedge}(K/p,f))=\left\{
 \begin{array}{rl}
0  & \text{if $n$ is odd }\\
(Z/(p^\infty))^k & \text{if $n$ is even}
\end{array}\right.$$

This ends the calculation for all odd primes. The homotopy groups of $\THH^{K_p^\wedge}(K/p)$ are 0 in odd degrees and $(Z/(p^\infty))^k$ in even degrees, where $k$ is a number between 1 and $p-1$ depending on the $A_\infty$ structure on $K/p$. This result was proved before by Angeltveit (\cite{A1}). He used the B\"okstedt spectral sequence to calculate Topological Hochschild Homology. 

\rem : This is the calculation identifying $K/p$ as the Thom spectrum of $S^1$. A similar calculation can be carried out for the Thom spectrum of $S^3$ to get the same results.


\end{document}